\documentclass[11pt]{article}

\usepackage{amsmath}
\usepackage{amssymb}
\usepackage{indentfirst}
\usepackage{graphics} 
\usepackage{color}
\makeatletter

\@addtoreset{equation}{section}
\makeatother
\def\<{\langle}
\def\>{\rangle}

\newtheorem{lem}{Lemma}[section]
\newtheorem{theo}{Theorem}[section]
\newtheorem{rem}{Remark}[section]
\newtheorem{pro}{Proposition}[section]

\makeatletter
   
   \@addtoreset{equation}{section}
\makeatother

\begin{document}
\title{\bf Fast energy decay for $2$-D wave equation with\\ localized damping near spatial infinity}
\author{Ryo Ikehata\thanks{Corresponding author: ikehatar@hiroshima-u.ac.jp} \\ {\small Department of Mathematics}\\ {\small Graduate School of Education} \\ {\small Hiroshima University} \\ {\small Higashi-Hiroshima 739-8524, Japan}}
\date{}
\maketitle
\begin{abstract}
We consider the Cauchy problem for wave equations with localized damping in ${\bf R}^{2}$. The damping is effective only near spatial infinity. We obtain fast energy decay estimate such that $O(t^{-2}\log t)$ as $t \to \infty$. Unlike the results for the two-dimensional exterior mixed problem case (\cite{Ike-0}), the difficulty of not being able to use Hardy-type inequalities is overcome by using Poincar\'e-type inequalities in all spaces and the finite propagation property of the solution to construct an estimate formula. In the two-dimensional case, when comparing the problem in the whole space with that in the exterior domain, we find that there is a significant difference in the sense that the former requires a logarithmic correction to the energy decay rate.
\end{abstract}
\section{Introduction}
\footnote[0]{Keywords and Phrases: Wave equation; $2$-D space; localized damping; total energy; fast decay, finite propagation property, multiplier method.}
\footnote[0]{2010 Mathematics Subject Classification. Primary 35L05; Secondary 35B45, 35L15.}
We consider the Cauchy problem associated to the wave equation with variable coefficient  in ${\bf R}^{2}$ as follows:
\begin{equation}
u_{tt}(t,x) -\Delta u(t,x) + a(x)u_{t}(t,x) = 0,\ \ \ (t,x)\in (0,\infty)\times {\bf R}^{2},\label{eqn}
\end{equation}
\begin{equation}
u(0,x)= u_{0}(x),\ \ u_{t}(0,x)= u_{1}(x),\ \ \ x\in{\bf R}^{2} ,\label{initial}
\end{equation}
where $(u_{0},u_{1})$ are initial data chosen as (for simplicity) $u_{0} \in C_{0}^{\infty}({\bf R}^{2})$, $u_{1} \in C_{0}^{\infty}({\bf R}^{2})$, and for some $R > 0$
\begin{equation}\label{w1}
{\rm supp}\,u_{j} \subset B_{R}(0)\quad (j = 0,1)
\end{equation}
with $B_{R}(0) := \{x \in {\bf R}^{2}\,:\,\vert x \vert \leq R\}$. Furthermore, the given function $a: {\bf R}^{2} \to {\bf R}$ satisfies the two assumptions below:\\
{\bf (A-1)}\, $a \in {\rm C}^{2}({\bf R}^{2})$, $a(x) \geq 0$ ($\forall x \in {\bf R}^{2}$), and $\vert a(x)\vert$ is bounded in ${\bf R}^{2}$,\\
{\bf (A-2)}\,there exist constants $L > 0$ and $\varepsilon_{0} > 0$ such that $a(x) \geq \varepsilon_{0}$ for $x \in {\bf R}^{2}$ satisfying $\vert x\vert > L$.

Here, we have set
\[u_{t}=\frac{\partial u}{\partial t},\quad  u_{tt}=\frac{\partial^2 u}{\partial t^2}, \quad \Delta = \sum_{j=1}^{2}\frac{\partial^{2}}{\partial x_{j}^{2}}, \quad x = (x_1,x_2).\]

Note that solutions and/or functions considered in this paper are all real valued.  

Considering the previous assumptions on the initial data and  $a(x)$ it is known that the problem \eqref{eqn}-\eqref{initial} has a unique smooth (classical) solution $u(t,x)$, which satisfies the finite propagation property (cf. \cite{Ikawa}):
\begin{equation}\label{w2}
u(t,x) = 0 \quad \textstyle{{\rm for}} \quad \vert x\vert > R + t\quad (t \geq 0).
\end{equation}

It is well-known that under the condition {\rm (A-1)} the following energy identity holds.
\begin{equation}\label{00}
E_{u}(t) + \int_{0}^{t}\int_{{\bf R}^{2}}a(x)\vert u_{s}(s,x)\vert^{2}dxds = E_{u}(0),
\end{equation}
where the total energy $E_{u}(t)$ to the equation (1.1) is defined by
\[E_{u}(t) :=\frac{1}{2}\int_{{\bf R}^{2}}\big(\vert u_{t}(t,x)\vert^{2} + \vert\nabla u(t,x)\vert^{2}\big)dx.\]

Two-dimensional problems in unbounded domains for several evolution equations are interesting. Let us summarize this topic through the two-dimensional wave equation. In particular, we will discuss the decay rate of the $L^2$-norm of solution itself and total energy.

First, we should cite Matsumura's celebrated result \cite{Ma} regarding energy decay in the following equation:
\[u_{tt}-\Delta u + u_{t} = 0.\]
The author \cite{Ma} reported the following optimal estimates in the two dimensional case.
\begin{equation}\label{L1}
\Vert u(t,\cdot)\Vert^{2} = O(t^{-1}),\quad E_{u}(t) = O(t^{-2})\,\, (t \to \infty),
\end{equation}
where $\Vert v\Vert$ denotes the usual $L^{2}$-norm of $v \in L^{2}({\bf R}^{2})$. This type of estimate for linear equation is also important in applications to nonlinear problems (see \cite{TY}), and as a linear problem in itself, it has become a target for subsequent research in the case when the damping $a(x)$ is sufficiently effective. 

The next issue to consider is the generalization of Matsumura's results to cases with variable friction coefficients $a(x)$ such as our equation \eqref{eqn}. Depending on the effectiveness of friction, the uniform decay of total energy becomes a delicate issue, and it is also interesting because it depends on the geometric shape of the (exterior) domain. A typical example of such a problem is the initial-boundary value problem of \eqref{eqn} in the smooth exterior domain, which has been studied extensively in previous research. For these reserch articles, one can cite \cite{R, D, M, MN, N, IM, Ike-0, LSM, DA} and the references therein. In this context, paper \cite{MN} (see also \cite{M}) in particular examines the exterior mixed problem of \eqref{eqn} outside the star-shaped region and derives the energy decay order of $O(t^{-1})$ at most for spatial dimension $n \geq 3$. On the other hand, in paper \cite{N} the author removed the geometric restriction of the star-shaped outer region using the so-called Lions condition and derived that the total energy decays at the order of $O(t^{-1})$ in cases of two or more spatial dimensions. When comparing Matsumura's two-dimensional results with the results of these exterior problems, it is natural to wonder whether there is any difference in the rate of energy decay between exterior problems and initial value problems in the entire space. Although this is a limited result in the star-shaped exterior domain, in paper \cite{Ike-0} the author revised partially Nakao's estimate using an improved version of Morawetz's method (see \cite{Mora}) and derived the fast decay rate of total energy $O(t^{-2})$ and $L^2$-norm in order $O(t^{-1})$. These results are being studied under conditions {\rm (A-1)} and {\rm (A-2)} in this paper. Although this is a result for star-shaped outer regions, it derives an estimate that agrees with Matsumura's $2$-Dimensional result \cite{Ma} even in the case of variable coefficients. Incidentally, in the case of exterior mixied problems with constant coefficients $a(x) = 1$ in \eqref{eqn}, it should be noted that Paper \cite{IM} obtained faster energy decay like $O(t^{-2})$, which was completely consistent with the local energy decay rate in the two-dimensional exterior mixied problem in Paper \cite{D}. Incidentaly, for research on diffusion phenomena in cases of friction with asymptotically periodic variable coefficients or friction that disappears at an appropriate order in space, we refer to the research of \cite{JJ} and \cite{SW} and the references cited therein.

Subsequently, papers \cite{LSM} and \cite{DA} were published, which used the so-called GCC condition on the damper $a(x)$ to generalize the exterior problem of the star-shaped region that had been restricted until then, and studied optimal rates of decay for the local and  total energy, and the $L^{2}$-norm of the solution. In particular, the results of \cite{DA} generalize those of \cite{Ike-0} by using GCC condition on the damper. From these results, one can roughly conclude that in the case of a two-dimensional exterior mixed problem, if friction is effective in the far field, and has some additional assumptions on the damper it is possible to reach the Matsumura estimate \eqref{L1}. The fact that the Hardy inequality holds in two dimensionnal exterior domain is important in some cases (cf. \cite{D, Ike-0}).

In cases where friction is effective only in distant spaces, the validity of Matsumura's estimate of the initial value problem \eqref{eqn}-\eqref{initial} in the entire space is a delicate issue. This is because it is difficult to control the $L^2$-norm of the solution near the origin. This is pointed our recently in \cite{Ike-1}, in fact, in \cite{Ike-1} the author derived for the free waves problem
\begin{equation}
v_{tt}(t,x) -\Delta v(t,x) = 0,\ \ \ (t,x)\in (0,\infty)\times {\bf R}^{2},\label{feqn}
\end{equation}
\begin{equation}
v(0,x)= u_{0}(x),\ \ v_{t}(0,x)= u_{1}(x),\ \ \ x\in{\bf R}^{2} ,\label{finitial}
\end{equation}
that the solution $v(t,x)$ to problem \eqref{feqn}-\eqref{finitial} satisfies
\[\Vert v(t,\cdot)\Vert \sim \sqrt{\log t} \quad (t \to \infty)\]
provided that $\int_{{\bf R}^{2}}u_{1}(x)dx \ne 0$. That is, the $L^{2}$-norm of the solution grows up to infinity with $\sqrt{\log t}$-order as time goes to infinity. Therefore, if the effect of friction $a(x)$ is poor near the origin, the singularity caused by such low frequencies in the solution itself will inevitably appear. This will be a major difference between exterior problems and all-space problems. 

The main objective of this paper is to capture the influence of free waves when friction $a(x)$ is weak near the origin. In particular, in the initial value problem of two-dimensional waves with friction $a(x)$ localized only in the far field, the question of whether the total energy decay rate is the same as Matsumura's estimate has been left unanswered for the past $20$ years. The ability to interpret the logarithmic growth order through prior research \cite{Ike-1} has provided valuable motivation for summarizing this paper.

Our main result is to report the following crucial result.
\begin{theo}\label{th1}\,Assume \eqref{w1}, {\rm (A-1)} and {\rm (A-2)}. Then the unique smooth solution $u(t,x)$ to problem {\rm (1.1)}-{\rm (1.2)} satisfies 
\[\int_{{\bf R}^{2}}\vert u(t,x)\vert^{2}dx = O(t^{-1}\log t)\quad (t \to \infty),\]
\[E_{u}(t) = O(t^{-2}\log t)\quad (t \to \infty).\]
\end{theo}
\par
\vspace{0.2cm}
\begin{rem}{\rm About the $L^{2}$-decay result as $t \to \infty$, near the origin of space (where friction has little effect), the influence of free waves discussed in \cite{Ike-1} can be observed (logarithmic growth in $t$ estimate appears), and in distant space, the influence of damped waves derived in \cite{Ma} appears ($t^{-1}$-decay appears), and the interaction between them supports the theorem's assertion.}
\end{rem}
\begin{rem}{\rm The obtained results in Theorem \ref{th1} completely improve the decay rates previously obtained in \cite[(1.2) and (1.3) of Theorem 1]{M} partially in the $2$-D case.}
\end{rem}
\begin{rem}{\rm The next problem is to remove the compact support condition imposed on the initial data. This may be a little difficult in the case of two-dimensional whole space.}
\end{rem}
This paper is organized as follows. In Section 2 we prepare important several lemmas, which will be used later, and we shall prove Theorem \ref{th1} in Section 3.\\

{\bf Notation.} {\small Throughout this paper, $\| \cdot\|_q$ stands for the usual $L^q({\bf R}^{n})$-norm. For simplicity of notation, in particular, we use $\| \cdot\|$ instead of $\| \cdot\|_2$. Furthermore, we denote $\Vert\cdot\Vert_{H^{1}}$ as the usual $H^{1}$-norm. 
The inner product in $L^{2}({\bf R}^{2})$ is denoted by
\[(f,g) := \int_{{\bf R}^{n}}f(x)g(x)dx, \quad f,g \in L^{2}({\bf R}^{n}).\] 
Moreover, we set
$$B_{r}(p) := \{x \in {\bf R}^{2}\,:\,\vert x - p \vert \leq r\}$$
for $p \in {\bf R}^{2}$ and $r > 0$.}

\section{$L^2$-bounds of solutions}

We will prove the main theorem in Section 3. In doing so, we will provide some supplementary lemmas in this Section 2.\\
First we prepare the so-called Poincar\'e inequality in the whole space (see \cite[the proof of Lemma 2.1 at page 352]{M}).
\begin{lem}\label{w10} For each $\rho > 0$ there exists a constant $C > 0$ depending only on $\rho$ such that
\[\int_{\vert x\vert \leq \rho}\vert v(x)\vert^{2}dx \leq C\left(\int_{{\bf R}^{2}}\vert\nabla v(x)\vert^{2}dx + \int_{\vert x\vert \geq \rho}\vert v(x)\vert^{2}dx\right)\]
for all $v \in H^{1}({\bf R}^{2})$. 
\end{lem}
Next, the following $L^{2}$-estimate can be used to prove our main theorem. The essential idea of the following lemma has already been developed in the analysis of the $2$-D elastic wave equation case in the preceding paper \cite{RR}, but a proof in the form of a reworked version of the equation \eqref{eqn} in this paper is included with a slight modification for the reader's convenience.
\begin{lem}\label{theorem7} Suppose {\rm (1.3)} and {\rm (A-1)}. The smooth solution $u(t,x)$ to problem \eqref{eqn}-\eqref{initial} satisfies the following growth property.
\[ \Vert u(t,\cdot)\Vert^{2} + \int_{0}^{t}\int_{{\bf R}^{2}}a(x)\vert u(s,x)\vert^{2}dxds \leq I_{0}(u_{0},u_{1},t)\quad (t \gg 1),\]
where 
\[I_{0}(u_{0},u_{1},t) := \Vert u_{0}\Vert^{2} + C_{1}\Vert u_{1} + a(\cdot)u_{0}\Vert_{1}^{2}\log t + C_{2}\Vert u_{1} + a(\cdot)u_{0}\Vert_{q}^{2},\]
and $C_{j} > 0$ {\rm (}$j = 1,2${\rm )} are some constants depending on $R$, and $q \in (2,\infty]$ is an arbirary fixed real number.
\end{lem}
\begin{rem}{\rm Since $\Vert u(t,\cdot)\Vert \leq C\sqrt{\log t}$ for $t \gg 1$ holds for the $L^2$-norm of the solution $u(t,x)$ to prolem \eqref{eqn}-\eqref{initial}, the $L^2$ estimate of this lemma seems to be the optimal when compared with the results of free waves obtained in \cite[Theorem 1.2]{Ike-1}.}

\end{rem}

Now, let us prove Lemma \ref{theorem7}. \\

For the smooth solution $u(t,x)$ to problem \eqref{eqn}-\eqref{initial}, as in \cite{IM} we set 
$$v(t,x) := \int_{0}^{t}u(s,x)ds.$$
Then, the function $v(t,x)$ satisfies
\begin{align}
& v_{tt}(t,x) - \Delta v(t,x) + a(x)v_{t}(t,x) = u_{1}(x) + a(x)u_{0}(x),\ \ \ (t,x)\in (0,\infty)\times {\bf R}^{2},\label{6eqn}\\
& v(0,x)= 0, \quad  v_{t}(0,x)= u_{0}(x),\ x\in{\bf R}^{2},\label{6initial}
\end{align}
and 
\[{\rm supp}\,v(t,\cdot) \subset B_{R+t}(0)\quad (\forall t \geq 0).\]
Now, we set 
\begin{equation}\label{N1}
h(x) := -\frac{1}{2\pi}\int_{{\bf R}^{2}}\log(\vert x- y\vert)(u_{1}(y)+a(y)u_{0}(y))dy,
\end{equation}
which represents the two dimensional Newton ptential. Since $u_{1}+a(\cdot)u_{0} \in C_{0}^{2}({\bf R}^{2})$, it follows from \cite[page 23, Theorem 1]{E} that $h \in C^{2}({\bf R}^{2})$, and the function $h(x)$ satisfies the Poisson equation:
\begin{equation}\label{N2}
-\Delta h(x) = u_{1}(x) + a(x)u_{0}(x), \quad x \in {\bf R}^{2}. 
\end{equation}  

We first prepare the following lemma already checked in \cite[Lemma 2.5]{RR}. Although this is a well-known result, it is recorded here for the convenience of readers.
\begin{lem}\label{N3}
The function $h(x)$ defined in \eqref{N1} satisfies
\[\vert x\vert \vert \nabla h(x)\vert \leq C\Vert u_{1} + a(x)u_{0}\Vert_{1}\]
for $\vert x \vert \geq 2R$, where $C > 0$ is an universal constant.
\end{lem}
{\it Proof of Lemma \ref{N3}.}\,\,First of all, note that under the assumption $\vert x\vert \geq 2R$ we can get
\[\vert x-y\vert \geq \vert x\vert - R \geq R,\quad \vert x\vert -R \geq \frac{1}{2}\vert x\vert\]
for $y \in {\bf R}^{2}$ satisfying $\vert y\vert \leq R$. After simple elementary computations on \eqref{N1} and the  assumption \eqref{w1} one can get the estimate
\[\vert\nabla h(x)\vert \leq \frac{1}{2\pi}\int_{{\bf R}^{2}}\frac{\vert u_{1}(y)+ a(y)u_{0}(y)\vert}{\vert x-y\vert}dy\]
\[= \frac{1}{2\pi}\int_{\vert y\vert \leq R}\frac{\vert u_{1}(y) + a(y)u_{0}(y)\vert}{\vert x-y\vert}dy \leq \frac{1}{2\pi}\frac{1}{\vert x\vert-R}\int_{{\bf R}^{2}}\vert u_{1}(y) + a(y)u_{0}(y)\vert dy\]
\[\leq \frac{1}{\pi}\frac{1}{\vert x\vert}\Vert u_{1} + a(\cdot)u_{0}\Vert_{1}.\]
This implies the desired estimate with a constant $C := \frac{1}{\pi}$.
\hfill
$\Box$\\

Let us prove Lemma \ref{theorem7} basing on Lemma \ref{N3}.\\
{\it Proof of Lemma \ref{theorem7}.}\,Now, multiplying both sides of \eqref{6eqn} by $v_{t}$, and integrating over $[0,t]\times {\bf R}^{2}$ it follows that
\[\frac{1}{2}\Vert v_{t}(t,\cdot)\Vert^{2} + \frac{1}{2}\Vert\nabla v(t,\cdot)\Vert^{2} + \int_{0}^{t}\int_{{\bf R}^{2}}a(x)\vert v_{s}(s,x)\vert^{2}dxds\]
\begin{equation}\label{N4}
= \frac{1}{2}\Vert u_{0}\Vert^{2} + \int_{{\bf R}^{2}}(u_{1}(x)+a(x)u_{0}(x))v(t,x)dx.
\end{equation}

Then, it follows from the finite speed of propagation property for $v(t,x)$, \eqref{N2} and the integration by parts that
\[\left\vert \int_{{\bf R}^{2}}(u_{1}(x)+a(x)u_{0}(x))v(t,x)dx\right\vert = \left\vert \int_{\vert x\vert \leq 2R + t}(u_{1}(x)+a(x)u_{0}(x))v(t,x)dx \right\vert\]
\[= \left\vert -\int_{\vert x\vert \leq 2R + t}\Delta h(x)v(t,x)dx \right\vert = \left\vert \int_{\vert x\vert \leq 2R + t}\nabla h(x)\cdot\nabla v(t,x)dx \right\vert\]
\[\leq \int_{\vert x\vert \leq 2R + t}\vert\nabla h(x)\vert\vert\nabla v(t,x)\vert dx\]
\[\leq \frac{1}{2}\int_{\vert x\vert \leq 2R + t}\vert\nabla h(x)\vert^{2}dx + \frac{1}{2}\int_{\vert x\vert \leq 2R + t}\vert\nabla v(t,x)\vert^{2}dx\]
\begin{equation}\label{N5}
\leq \frac{1}{2}\int_{\vert x\vert \leq 2R + t}\vert\nabla h(x)\vert^{2}dx + \frac{1}{2}\int_{{\bf R}^{2}}\vert\nabla v(t,x)\vert^{2}dx.
\end{equation} 
To get \eqref{N5} one has just used the facts that $v(t,x) = 0$, and $\vert\nabla h(x)\vert < +\infty$ for $\vert x\vert = 2R + t$ to the boundary integral.
Thus, \eqref{N4} and \eqref{N5} imply
\[\frac{1}{2}\Vert v_{t}(t,\cdot)\Vert^{2} + \int_{0}^{t}\int_{{\bf R}^{2}}a(x)\vert v_{s}(s,x)\vert^{2}dxds\]
\begin{equation}\label{N6}
\leq \frac{1}{2}\Vert u_{0}\Vert^{2} + \frac{1}{2}\int_{\vert x\vert \leq 2R + t}\vert\nabla h(x)\vert^{2}dx. 
\end{equation}
Note that 
\[\int_{\vert x\vert \leq 2R + t}\vert\nabla h(x)\vert^{2}dx < +\infty\]
for each $t \geq 0$ because of $h \in C^{2}({\bf R}^{2})$. Therefore, one can get the crucial estimate because of $v_{t} = u$.
\begin{equation}\label{N7}
\Vert u(t,\cdot)\Vert^{2} + \int_{0}^{t}\int_{{\bf R}^{2}}a(x)\vert u(s,x)\vert^{2}dxds \leq \Vert u_{0}\Vert^{2} + \int_{\vert x\vert \leq 2R + t}\vert\nabla h(x)\vert^{2}dx.
\end{equation}

Let us estimate the second term of the right hand side of \eqref{N7} by using Lemma \ref{N3}. We set
\[I_{h} := \int_{\vert x\vert \leq 2R}\vert\nabla h(x)\vert^{2}dx.\]
We first note the following decomposition. 
\begin{equation}\label{N8}
\int_{\vert x\vert \leq 2R + t}\vert\nabla h(x)\vert^{2}dx = I_{h} + \int_{2R \leq \vert x\vert \leq 2R + t}\vert\nabla h(x)\vert^{2}dx.
\end{equation}
Then, from Lemma \ref{N3} it follows that
\[\int_{2R \leq \vert x\vert \leq 2R + t}\vert\nabla h(x)\vert^{2}dx \leq C^{2}\Vert u_{1} + a(\cdot)u_{1}\Vert_{1}^{2}\int_{2R \leq \vert x\vert \leq 2R + t}\frac{1}{\vert x\vert^{2}}dx\]
\[= 2\pi C^{2}\Vert u_{1} + a(\cdot)u_{1}\Vert_{1}^{2}\int_{2R}^{2R + t}\frac{r}{r^{2}}dr\]
\[= 2\pi C^{2}\Vert u_{1} + a(\cdot)u_{1}\Vert_{1}^{2}\left(\log(2R+t)-\log(2R)\right)\]
\begin{equation}\label{N9}
\leq 2\pi C^{2}\Vert u_{1} + a(\cdot)u_{1}\Vert_{1}^{2}\log(2R+t)\quad (t \geq 0). 
\end{equation}
Therefore, by \eqref{N7}, \eqref{N8} and \eqref{N9} one has arrived at the crucial estimate:
\begin{equation}\label{N10}
\Vert u(t,\cdot)\Vert^{2} + \int_{0}^{t}\int_{{\bf R}^{2}}a(x)\vert u(s,x)\vert^{2}dxds \leq \Vert u_{0}\Vert^{2} + I_{h} + 2\pi C^{2}\Vert u_{1}+a(\cdot)u_{0}\Vert_{1}^{2}\log(2R+t) \quad (t \geq 0).
\end{equation}
\noindent
Let us estimate $I_{h}$ in terms of $u_{1}(x) + a(x)u_{1}(x)$. Note that generally it holds that
\[\int_{{\bf R}^{2}}\vert\nabla h(x)\vert^{2}dx = +\infty.\] 

For an arbitrarily fixed $x_{0} \in B_{2R}(0)$, we can easily arrive at the estimate basing on the H\"older inequality.
\[\vert\nabla h(x_{0})\vert \leq \frac{1}{2\pi}\int_{{\bf R}^{2}}\frac{\vert u_{1}(y) + a(y)u_{0}(y)\vert}{\vert x_{0}-y\vert}dy = \frac{1}{2\pi}\int_{\vert y\vert \leq R}\frac{\vert u_{1}(y)+ a(y)u_{0}(y)\vert}{\vert x_{0}-y\vert}dy \]
\[= \frac{1}{2\pi}\int_{B_{R}(x_{0})}\frac{\vert u_{1}(x_{0}-z)+ a(x_{0}-z)u_{0}(x_{0}-z)\vert}{\vert z\vert}dz\]
\begin{equation}\label{N11}
\leq \frac{1}{2\pi}\left(\int_{B_{R}(x_{0})}\frac{1}{\vert z\vert^{p}}dz\right)^{\frac{1}{p}}\Vert u_{1} + a(\cdot)u_{0}\Vert_{q},
\end{equation}
where $1/p + 1/q = 1$, and $p,q \in [1,\infty]$ at this stage. One notes that in the case of $x_{0} \in B_{2R}(0)$, one has, for example, $B_{R}(x_{0}) \subset B_{4R}(0)$. Therefore, we get
\begin{equation}\label{N12}
\int_{B_{R}(x_{0})}\frac{1}{\vert z\vert^{p}}dz \leq \int_{B_{4R}(0)}\frac{1}{\vert z\vert^{p}}dz = \frac{2\pi}{2-p}(4R)^{2-p},
\end{equation}
provided that $p \in [1,2)$. Put $x_{0}$ back into $x$ to get the following by \eqref{N11} and \eqref{N12}. 
\[\vert\nabla h(x)\vert \leq (2\pi)^{\frac{1}{p}-1}(2-p)^{-\frac{1}{p}}(4R)^{\frac{2-p}{p}}\Vert u_{1}+a(\cdot)u_{0}\Vert_{q}\quad \forall x \in B_{2R}(0).\]
Therefore, one can obtain the desired estimate:
\begin{equation}\label{N13}
I_{h} = \int_{\vert x\vert \leq 2R}\vert\nabla h(x)\vert^{2}dx \leq C_{R}\Vert u_{1}+a(\cdot)u_{0}\Vert_{q}^{2}
\end{equation}
with some constant $C_{R} > 0$ defined by
\[C_{R} := 4\pi R^{2}\{(2\pi)^{\frac{1}{p}-1}(2-p)^{-\frac{1}{p}}(4R)^{\frac{2-p}{p}}\}^{2}.\]
Here, one should impose the assumption on $p$ and $q$ such that
\[p \in [1,2),\quad q := \frac{p}{p-1} \in (2,\infty].\] 
Of course, from the regularity assumed on the initial data we see that $\Vert u_{1}+a(\cdot)u_{0}\Vert_{q} < +\infty$ for any $q \in (2,\infty]$.

At last one can arrive at the desired estimate in Lemma \ref{theorem7} because of \eqref{N10} and \eqref{N13}
\hfill
$\Box$

\section{Proof of Theorem \ref{th1}}

Under preparations as in Section 2, let us prove Theorem \ref{th1}. The strategy is to employ the multiplier methods developed in \cite{Z}, \cite{N} and \cite{Ike-0}, and use the preliminary estimates in Section 2 where necessary. In particular, we borrow Nakao's computations in the first part of computations.\\
Set
\[ \phi(r) = \left\{
    \begin{array}{ll}
     \displaystyle{\varepsilon_{0}}&
           \qquad 0 \leq r := \vert x\vert \leq L, \\[0.2cm]
    \displaystyle{\frac{\varepsilon_{0}L}{r}}& \qquad r \geq L.
    \end{array} \right. \]
With a smooth vector field $h(x) := (h_{1}(x),h_{2}(x))$, we start with the following identity.
\[\frac{d}{dt}\int_{{\bf R}^{2}}u_{t}(t,x)\left(h(x)\cdot\nabla u(t,x)\right)dx + \frac{1}{2}\int_{{\bf R}^{2}}\nabla\cdot h(x)\left(\vert u_{t}(t,x)\vert^{2}-\vert\nabla u(t,x)\vert^{2}  \right)dx\] 
\begin{equation}\label{M1}
+ \sum_{j,k = 1}^{2}\int_{{\bf R}^{2}}u_{x_{j}}u_{x_{k}}\frac{\partial}{\partial x_{j}}h_{k}(x)dx + \int_{{\bf R}^{2}}a(x)u_{t}(t,x)\left(h(x)\cdot \nabla u(t,x)\right)dx = 0.
\end{equation} 
By substituting $h(x) := \phi(r)x$ in \eqref{M1} we see that
\[\frac{d}{dt}\int_{{\bf R}^{2}}u_{t}(t,x)\phi(r)(x\cdot\nabla u(t,x))dx + \frac{1}{2}\int_{{\bf R}^{2}}(2\phi(r)+\phi'(r)r)\left(\vert u_{t}(t,x)\vert^{2}-\vert\nabla u(t,x)\vert^{2}  \right)dx\] 
\begin{equation}\label{M2}
+ \int_{{\bf R}^{2}}\left(\frac{\phi'(r)}{r}\vert x\cdot\nabla u(t,x)\vert^{2} + \phi(r)\vert\nabla u(t,x)\vert^{2}\right)dx + \int_{{\bf R}^{2}}a(x)u_{t}(t,x)\phi(r)(x\cdot \nabla u(t,x)) dx = 0.
\end{equation}
On the other hand, we have the following identity for $\alpha > 0$.
\begin{equation}\label{M2-1}
\alpha\frac{d}{dt}\left(u_{t}(t,\cdot),u(t,\cdot) \right) - \alpha\Vert u_{t}(t,\cdot)\Vert^{2} + \alpha\Vert \nabla u(t,\cdot)\Vert^{2} + \frac{\alpha}{2}\frac{d}{dt}\int_{{\bf R}^{2}}a(x)\vert u(t,x)\vert^{2}dx = 0.
\end{equation}
\eqref{M2-1} can be easily derived by multiplying both sides of \eqref{eqn} by $\alpha u(t,x)$, and integrating over ${\bf R}^{2}$. Therefore, adding three identities \eqref{M2}, \eqref{M2-1} and (\eqref{00} $\times k$) ($k > 0$) yield the following identity.
\[\frac{d}{dt}\left( \int_{{\bf R}^{2}}u_{t}(t,x)\phi(r)(x\cdot \nabla u(t,x))dx + \alpha(u_{t}(t,\cdot),u(t,\cdot)) + \frac{\alpha}{2}\int_{{\bf R}^{2}}a(x)\vert u(t,x)\vert^{2}dx + kE_{u}(t) \right)\]
\[+ \int_{{\bf R}^{2}}\left(\frac{2\phi(r)+\phi'(r)r}{2}-\alpha + ka(x) \right)\vert u_{t}(t,x)\vert^{2}dx\]
\[+ \int_{{\bf R}^{2}}\left( \alpha - \frac{2\phi(r)+\phi'(r)r}{2} + \phi(r) + \phi'(r)r \right)\vert\nabla u(t,x)\vert^{2}dx\]
\begin{equation}\label{M3}
= - \int_{{\bf R}^{2}}a(x)u_{t}(t,x)\phi(r)(x\cdot \nabla u(t,x))dx,
\end{equation}
where the constants $\alpha > 0$ and $k > 0$ again can be choosen suitably later on. By assumptions (A-1) and (A-2), one can choose large $\alpha > 0$ such that for any $k \geq 3$
\begin{equation}\label{M4}
\frac{2\phi(r)+ \phi'(r)r}{2} + \frac{ka(x)}{2} - \alpha > \varepsilon_{1} > 0,
\end{equation} 
\begin{equation}\label{M5}
\alpha - \frac{2\phi(r)+ \phi'(r)r}{2} + \phi(r) + \phi'(r)r > \varepsilon_{1} > 0,
\end{equation} 
where $\varepsilon_{1} >0$ satisfies
\[0 < \varepsilon_{1} < \frac{\varepsilon_{0}}{2}.\]
In fact, as in \cite{N, Ike-0} one can choose
\[\alpha := \frac{3}{4}\varepsilon_{0}, \quad \varepsilon_{1} := \frac{\varepsilon_{0}}{8}\] 
in the two dimensional case. Thus, it holds that
\[\frac{d}{dt}\left( \int_{{\bf R}^{2}}u_{t}(t,x)\phi(r)(x\cdot \nabla u(t,x))dx + \alpha(u_{t}(t,\cdot),u(t,\cdot)) + \frac{\alpha}{2}\int_{{\bf R}^{2}}a(x)\vert u(t,x)\vert^{2}dx + kE_{u}(t) \right)\]
\[+ \frac{k}{2}\int_{{\bf R}^{2}}a(x)\vert u_{t}(t,x)\vert^{2}dx + \varepsilon_{1}\int_{{\bf R}^{2}}\vert\nabla u(t,x)\vert^{2}dx\]
\begin{equation}\label{M6}
\leq - \int_{{\bf R}^{2}}a(x)u_{t}(t,x)\phi(r)(x\cdot \nabla u(t,x))dx.
\end{equation}
While, with some constant $C_{3} > 0$ depending on $\Vert a\Vert_{\infty}$, $\varepsilon_{0}$ and $L$ it follows that
\begin{equation}\label{M7}
\vert \int_{{\bf R}^{2}}a(x)u_{t}(t,x)\phi(r)(x\cdot \nabla u(t,x))dx\vert \leq \frac{k}{4}\int_{{\bf R}^{2}}a(x)\vert u_{t}(t,x)\vert^{2}dx + \frac{C_{3}}{k}\int_{{\bf R}^{2}}\vert\nabla u(t,x)\vert^{2}dx
\end{equation}
for any $k > 0$. Therefore, from \eqref{M6} and \eqref{M7} it follows that
\[\frac{d}{dt}G_{k}(t) + \frac{k}{4}\int_{{\bf R}^{2}}a(x)\vert u_{t}(t,x)\vert^{2}dx + \varepsilon_{1}\int_{{\bf R}^{2}}(\vert u_{t}(t,x)\vert^{2} + \vert\nabla u(t,x)\vert^{2})dx \leq \frac{2C_{3}}{k}E_{u}(t),\]
where
\[G_{k}(t) := \int_{{\bf R}^{2}}u_{t}(t,x)\phi(r)(x\cdot \nabla u(t,x))dx + \alpha(u_{t}(t,\cdot),u(t,\cdot)) + \frac{\alpha}{2}\int_{{\bf R}^{2}}a(x)\vert u(t,x)\vert^{2}dx + kE_{u}(t),\]
so that by taking $k \gg 1$ large enough such that
\[\frac{C_{3}}{k} < \varepsilon_{1},\]
and setting $\beta := 2(\varepsilon_{1}-\displaystyle{\frac{C_{3}}{k}}) > 0$, one can arrive at the inequality.
\begin{equation}\label{M70}
\frac{d}{dt}G_{k}(t) + \beta E_{u}(t) \leq 0,
\end{equation}
which implies the important estimate.
\begin{equation}\label{M8}
G_{k}(t) + \beta\int_{0}^{t} E_{u}(s)ds \leq G_{k}(0) \leq C_{4}(\Vert u_{0}\Vert_{H^{1}}^{2} + \Vert u_{1}\Vert^{2}),
\end{equation}
where $C_{4} > 0$ is a constant depending on $\Vert a\Vert_{\infty}$, $\varepsilon_{0} > 0$ and $L > 0$. 

Let us present the following lemma to get the energy decay estimate. In order to prove Lemma \ref{M9} below one uses Lemma \ref{w10} stated in Section 2. 
\begin{lem}\label{M9}For large $k > 3$, it holds that
\[G_{k}(t) \geq 0 \quad (\forall t \geq 0).\]
\end{lem}
{\it Proof of Lemma \ref{M9}.}\, It follows from Lemma \ref{w10} and assumption {\rm (A-2)} that for an arbitrary $\varepsilon > 0$
\[-\alpha\left(u_{t}(t,\cdot),u(t,\cdot)\right) \leq \frac{\alpha}{2\varepsilon}\Vert u_{t}(t,\cdot)\Vert^{2} + \frac{\varepsilon\alpha}{2}\Vert u(t,\cdot)\Vert^{2}\] 
\[\leq \frac{\alpha}{\varepsilon}E_{u}(t) + \frac{\varepsilon\alpha}{2}\left( \int_{\vert x\vert \geq L}\vert u(t,x)\vert^{2}dx +  \int_{\vert x\vert \leq L}\vert u(t,x)\vert^{2}dx    \right)\]
\[\leq \frac{\alpha}{\varepsilon}E_{u}(t) + \frac{\varepsilon\alpha}{2\varepsilon_{0}}\left( \int_{\vert x\vert \geq L}a(x)\vert u(t,x)\vert^{2}dx \right) +  \frac{\varepsilon\alpha}{2}\int_{\vert x\vert \leq L}\vert u(t,x)\vert^{2}dx\]
\[\leq \frac{\alpha}{\varepsilon}E_{u}(t) + \frac{\varepsilon\alpha}{2\varepsilon_{0}}\int_{\vert x\vert \geq L}a(x)\vert u(t,x)\vert^{2}dx + \frac{\varepsilon\alpha C}{2}\left(\int_{{\bf R}^{2}}\vert\nabla u(t,x)\vert^{2}dx + \int_{\vert x\vert \geq L}\vert u(t,x)\vert^{2}dx\right)\]
\[\leq \frac{\alpha}{\varepsilon}E_{u}(t) + \frac{\varepsilon\alpha}{2\varepsilon_{0}}\int_{{\bf R}^{2}}a(x)\vert u(t,x)\vert^{2}dx +  \varepsilon\alpha CE_{u}(t) + \frac{\varepsilon\alpha C}{2\varepsilon_{0}}\int_{{\bf R}^{2}}a(x)\vert u(t,x)\vert^{2}dx\]
\begin{equation}\label{M10}
= (\frac{\alpha}{\varepsilon}+\varepsilon\alpha C)E_{u}(t) + (1+C)\frac{\varepsilon\alpha}{2\varepsilon_{0}}\int_{{\bf R}^{2}}a(x)\vert u(t,x)\vert^{2}dx.
\end{equation}
While, from the same estimate derived in \cite[page 397, proof of Lemma 2.3]{Ike-0} we can review that
\begin{equation}\label{M11}
-\int_{{\bf R}^{2}}u_{t}(t,x)\phi(r)(x\cdot\nabla u(t,x))dx \leq 2\varepsilon_{0}LE_{u}(t).
\end{equation}
By adding the two inequalities \eqref{M10} and \eqref{M11} one can arrive at the inequality.
\[-\alpha\left(u_{t}(t,\cdot),u(t,\cdot)\right) - \int_{{\bf R}^{2}}u_{t}(t,x)\phi(r)(x\cdot\nabla u(t,x))dx \]
\begin{equation}\label{M12}
\leq (\frac{\alpha}{\varepsilon} + \varepsilon\alpha C + 2\varepsilon_{0}L)E_{u}(t) + (1+C)\frac{\varepsilon\alpha}{2\varepsilon_{0}}\int_{{\bf R}^{2}}a(x)\vert u(t,x)\vert^{2}dx.
\end{equation}
By choosing $\varepsilon> 0$ small satisfying $\varepsilon < \varepsilon_{0}/(1+C)$, and taking $k > 3$ large enough such that 
$$(\frac{\alpha}{\varepsilon} + \varepsilon\alpha C + 2\varepsilon_{0}L) < k,$$
one can get the desired estimate $G_{k}(t) \geq 0$ for $t \geq 0$. 
\hfill
$\Box$

It follows from \eqref{M8} and Lemma \ref{M9} one can easily derive the following proposition (see \cite[Lemma 2.4]{Ike-0}). This intermediate result shows the algebraic decay rate $O(t^{-1})$ ($t \to \infty$) of the total energy $E_{u}(t)$. 
\begin{pro}\label{M13}For large $k > 0$ it holds that
\[\int_{0}^{t}E_{u}(s)ds \leq \frac{1}{\beta}G_{k}(0)\quad (t \geq 0),\]
\[E_{u}(t) \leq \left(E_{u}(0) + \frac{1}{\beta}G_{k}(0)\right)(1+t)^{-1}\quad (t \geq 0).\]
\end{pro}

Next, by combining Lemma \ref{theorem7}, assumption {\rm (A-2)}, and Proposition \ref{M13} one can get the following estimate.
\begin{lem}\label{M14} It holds that for large $k > 0$
\[\int_{0}^{t}\vert\left(u_{s}(s,\cdot),u(s,\cdot)\right)\vert ds \leq \frac{1+C}{\beta}G_{k}(0) + \frac{1+C}{2\varepsilon_{0}}I_{0}(u_{0},u_{1},t) \quad (t \gg 1),\]
where $I_{0}(u_{0},u_{1},t) $ is a quantity already defined in Lemma {\rm \ref{theorem7}}.
\end{lem} 
{\it Proof.}\,We have the following series of estimates by using assumption {\rm (A-2)} and Lemma \ref{w10} in suitable places. 
\[\int_{0}^{t}\vert\left(u_{s}(s,\cdot),u(s,\cdot)\right)\vert ds\]
\[\leq \frac{1}{2}\int_{0}^{t}\Vert u_{s}(s,\cdot)\Vert^{2}ds + \frac{1}{2}\int_{0}^{t}\Vert u(s,\cdot)\Vert^{2}ds\]
\[\leq \int_{0}^{t}E_{u}(s)ds + \frac{1}{2}\int_{0}^{t}\int_{\vert x\vert \geq L}\vert u(s,x)\vert^{2}dxds + \frac{1}{2}\int_{0}^{t}\int_{\vert x\vert \leq L}\vert u(s,x)\vert^{2}dxds\]
\[\leq \int_{0}^{t}E_{u}(s)ds + \frac{1}{2\varepsilon_{0}}\int_{0}^{t}\int_{{\bf R}^{2}}a(x)\vert u(s,x)\vert^{2}dxds + \frac{C}{2}\int_{0}^{t}\left(\int_{{\bf R}^{2}}\vert \nabla u(s,x)\vert^{2}dx+\int_{\vert x\vert \geq L}\vert u(s,x)\vert^{2}dx\right)ds\]
\[\leq \int_{0}^{t}E_{u}(s)ds + \frac{1}{2\varepsilon_{0}}\int_{0}^{t}\int_{{\bf R}^{2}}a(x)\vert u(s,x)\vert^{2}dxds + C\int_{0}^{t}E_{u}(s)ds + \frac{C}{2\varepsilon_{0}}\int_{0}^{t}\int_{{\bf R}^{2}}a(x)\vert u(s,x)\vert^{2}dxds\]
\[=(1+C)\int_{0}^{t}E_{u}(s)ds + \frac{1+C}{2\varepsilon_{0}}\int_{0}^{t}\int_{{\bf R}^{2}}a(x)\vert u(s,x)\vert^{2}dxds.\]
The desired estimate is a direct consequence of  Lemma \ref{theorem7} and Proposition \ref{M13}.
\hfill
$\Box$

We also prepare the following lemma to finalize the proof of main Theorem \ref{th1}. 

\begin{lem}\label{M15} For large $k > 0$, it holds that\\
{\rm (i)}\,$\displaystyle{\int_{0}^{t}\int_{{\bf R}^{2}}}\left\vert u_{s}(s,x)\phi(r)(x\cdot\nabla u(s,x)) \right\vert dxds \leq 2\varepsilon_{0}L\beta^{-1}G_{k}(0)\quad (t \geq 0)$,\\ 
{\rm (ii)}\,$\displaystyle{\int_{0}^{t}}G_{k}(s)ds \leq K_{1}(t)$ for $t \gg 1$,\\
{\rm (iii)}\,$\displaystyle{\int_{0}^{t}}(1+s)E_{u}(s)ds \leq \beta^{-1}(G_{k}(0) + K_{1}(t))$ for $t \gg 1$, 
where
\[K_{1}(t) := \beta^{-1}(2\varepsilon_{0}L + \alpha(1+C) +k) + \frac{\alpha}{2}(1+\frac{1+C}{\varepsilon_{0}})I_{0}(u_{0},u_{1},t).\]
\end{lem}
{\it Proof.}\,Let's check them in order.\\
{\rm (i)} is a direct consequence of \cite[(2.13) at page 400]{Ike-0}.\\
{\rm (ii)}\,By definition of $G_{k}(t)$ we see that
\[\int_{0}^{t}G_{k}(s)ds \leq \int_{0}^{t}\int_{{\bf R}^{2}}\left\vert u_{s}(s,x)\phi(r)(x\cdot\nabla u(s,x)) \right\vert dxds + \alpha\int_{0}^{t}\vert\left(u_{s}(s,\cdot),u(s,\cdot)\right)\vert ds \]
\[+\frac{\alpha}{2}\int_{0}^{t}\int_{{\bf R}^{2}}a(x)\vert u(s,x)\vert^{2}dxds + k\int_{0}^{t}E_{u}(t)dt.\]
The desired estimate {\rm (ii)} can be checked by using {\rm (i)}, Lemma \ref{M14}, Lemma \ref{theorem7} and Proposition \ref{M13}. \\ 
{\rm (iii)} can be derived by multiplying both sides of  \eqref{M70} by $(1+t)$, then integrating over $[0,t]$ and further using {\rm (ii)} together with Lemma \ref{M9}, which is omitted.
\hfill
$\Box$

We will now begin the task of deriving the $L^2$ decay rate, which plays an essential role in achieving this rapid energy decay. We first evaluate the following quantity. 
$$V(t) := \int_{\vert x\vert \geq L}\vert u(t,x)\vert^{2}dx$$
in the exterior region in order to estimate  
$$U(t) := \int_{\vert x\vert \leq L}\vert u(t,x)\vert^{2}dx$$
in the interior region. Lemma \ref{w10} serves as a bridge between the two. The estimate of the interior domain comes down to the estimate of the exterior domain. The following is the growth estimate in the exterior region $\vert x\vert \geq L$. 
\begin{lem}\label{M15-1}For $k \gg 1$ and $t \gg 1$, it holds that
\[C_{5}(1+t)\int_{\vert x\vert \geq L}\vert u(t,x)\vert^{2}dx\]
\[\leq  K_{0} + \frac{1}{\varepsilon}(E_{u}(0)+ \frac{1}{\beta}G_{k}(0)) + \frac{2}{\beta}G_{k}(0) + (\frac{2}{\alpha}+\frac{2}{\beta})K_{1}(t) + \varepsilon C(E_{u}(0)+\frac{1}{\beta}G_{k}(0)) := K_{2}(t),\]
where $C_{5} > 0$ is an universal constant, and 
\[K_{0} := (u_{0},u_{1}) + \frac{1}{2}\int_{{\bf R}^{2}}a(x)\vert u_{0}(x)\vert^{2}dx.\] 
\end{lem}
{\it Proof.}\,We first note the following estimate based on Lemma \ref{w10} and Proposition \ref{M13}.
\[(1+t)U(t) \leq C(1+t)\int_{{\bf R}^{2}}\vert\nabla u(t,x)\vert^{2}dx + C(1+t)V(t)\]
\[\leq 2C(1+t)E_{u}(t) + C(1+t)V(t)\]
\begin{equation}\label{M16}
\leq 2C(E_{u}(0)+ \beta^{-1}G_{k}(0)) + C(1+t)V(t),
\end{equation}
for $t \geq 0$ and $k \gg 1$. Multiplying both sides of equation \eqref{eqn} by $(1+t)u(t,x)$ and integrating over $[0,t]\times {\bf R}^{2}$ in spacetime, we obtain the following.
\[\frac{1+t}{2}\int_{{\bf R}^{2}}a(x)\vert u(t,x)\vert^{2}dx + \frac{1}{2}\Vert u_{0}\Vert^{2} + \int_{0}^{t}(1+s)\Vert\nabla u(s,\cdot)\Vert^{2}ds\]
\begin{equation}\label{M17}
= K_{0} -(1+t)(u_{t}(t,\cdot),u(t,\cdot)) + \frac{1}{2}\Vert u(t,\cdot)\Vert^{2}+ \int_{0}^{t}(1+s)\Vert u_{s}(s,\cdot)\Vert^{2}ds + \frac{1}{2}\int_{0}^{t}\int_{{\bf R}^{2}}a(x)\vert u(s,x)\vert^{2}dxds,
\end{equation}
where
\[K_{0} := (u_{0},u_{1}) + \frac{1}{2}\int_{{\bf R}^{2}}a(x)\vert u_{0}(x)\vert^{2}dx.\]
Thus, by assumption {\rm (A-2)}, Lemma \ref{theorem7}, Proposition \ref{M13} and {\rm (iii)} of Lemma \ref{M15}, one has
\[\frac{1+t}{2}\varepsilon_{0}V(t) \leq K_{0} -(1+t)(u_{t}(t,\cdot),u(t,\cdot))\]
\[+ \frac{1}{2}\left(\Vert u(t,\cdot)\Vert^{2} + \int_{0}^{t}\int_{{\bf R}^{2}}a(x)\vert u(s,x)\vert^{2}dxds\right) + 2\int_{0}^{t}(1+s)E_{u}(s)ds\]
\[\leq K_{0} + \frac{(1+t)}{2\varepsilon}\Vert u_{t}(t,\cdot)\Vert^{2} + \frac{(1+t)\varepsilon}{2}\Vert u(t,\cdot)\Vert^{2}\]
\[+ I_{0}(u_{0},u_{1},t) + \frac{2}{\beta}G_{k}(0) + \frac{2}{\beta}K_{1}(t)\]
\[\leq K_{0} + \frac{1}{\varepsilon}(E_{u}(0)+ \frac{1}{\beta}G_{k}(0)) + \frac{2}{\beta}G_{k}(0) + \frac{(1+t)\varepsilon}{2}\Vert u(t,\cdot)\Vert^{2}\]
\[+\frac{2}{\beta}K_{1}(t) +I_{0}(u_{0},u_{1},t) \quad (t \gg 1),\]
where $\varepsilon > 0$ is a real number specified later on. Since
\[\Vert u(t,\cdot)\Vert^{2} = V(t) + U(t),\]
this inequality implies
\[\frac{(1+t)}{2}(\varepsilon_{0}-\varepsilon)V(t)\]
\[\leq K_{0} + \frac{1}{\varepsilon}(E_{u}(0)+ \frac{1}{\beta}G_{k}(0)) + \frac{2}{\beta}G_{k}(0) + \frac{2}{\beta}K_{1}(t) \]
\[+ I_{0}(u_{0},u_{1},t) + \frac{(1+t)\varepsilon}{2}U(t).\]
Applying \eqref{M16} to the end of the equation we just derived, one obtains the following.
\[\frac{(1+t)}{2}(\varepsilon_{0}-\varepsilon-C\varepsilon)V(t)\]
\[\leq K_{0} + \frac{1}{\varepsilon}(E_{u}(0)+ \frac{1}{\beta}G_{k}(0)) + \frac{2}{\beta}G_{k}(0) + \frac{2}{\beta}K_{1}(t) \]
\[+ I_{0}(u_{0},u_{1},t) + \varepsilon C(E_{u}(0)+\frac{1}{\beta}G_{k}(0)).\]
Since
\[I_{0}(u_{0},u_{1},t) \leq \frac{2}{\alpha}(1+\frac{1+C}{\varepsilon_{0}})^{-1}K_{1}(t),\]
by choosing $\varepsilon > 0$ small such that
\[C_{5} := (\varepsilon_{0}-\varepsilon-C\varepsilon)/2 > 0,\]
for such $\varepsilon > 0$ one has the desired estimate.
\hfill
$\Box$ 

Once Lemma \ref{M15-1} has been demonstrated, it can be easily derived the interior domain's estimate by using Proposition \ref{M13} and Lemma \ref{w10}. We state it without proof.
\begin{lem}\label{M18}\,For $k \gg 1$ and $t \gg 1$, it holds that
\[(1+t)U(t) \leq 2(E_{u}(0)+\beta^{-1}G_{k}(0)) + \frac{C}{C_{5}}K_{2}(t).\]
\end{lem}
By combining Lemmas \ref{M15-1} and \ref{M18} one has arrived at the following $L^{2}$-growth estimate of the solution to problem \eqref{eqn}-\eqref{initial}.
\begin{pro}\label{M19}\,Under the assumptions as in Theorem {\rm \ref{th1}}, for $k \gg 1$ and $t \gg 1$ it holds that
\[\Vert u(t,\cdot)\Vert^{2} = O(t^{-1}\log t)\quad (t \to \infty).\]
\end{pro}
{\it Proof.}\, It is trivial from the definition of the quantities $K_{1}(t)$ and/or $K_{2}(t)$, and the results derived in Lemmas \ref{M15-1} and \ref{M18}.
\hfill
$\Box$
\\
Finally, let us prove the total energy decay result.
\begin{pro}\label{M20}\,Under the assumptions as in Theorem {\rm \ref{th1}}, for $k \gg 1$ and $t \gg 1$ it holds that
\[E_{u}(t) = O(t^{-2}\log t)\quad (t \to \infty).\]
\end{pro}
{\it Proof.}\,We start with the following computation on the total energy.
\[\frac{d}{dt}\left((1+t)^{2}E_{u}(t)\right) = 2(1+t)E_{u}(t) + (1+t)^{2}E_{u}'(t)\] 
\[\leq 2(1+t)E_{u}(t)\quad (t \geq 0).\]
Integrating both sides of the inequality above over $[0,t]$ one can get the inequality:
\[(1+t)^{2}E_{u}(t) \leq E_{u}(0) + 2\int_{0}^{t}(1+s)E_{u}(s)ds.\]
Thus, by using {\rm (iii)} of Lemma \ref{M15} it holds that
\[(1+t)^{2}E_{u}(t) \leq E_{u}(0) + \frac{2}{\beta}G_{k}(0) + \frac{2}{\beta}K_{1}(t)\]
for large $k > 3$ and large $t > 0$. From the definition of $K_{1}(t)$, the desired result can be easily derived.
\hfill
$\Box$

Let us finalize the proof of Theorem \ref{th1}.\\
\noindent
{\it Proof of Theorem \ref{th1} completed.}\, The both decay estimates for $\Vert u(t,\cdot)\Vert^{2}$ and $E_{u}(t)$ has already been obtained by those of Propositions \ref{M19} and \ref{M20}.
\hfill
$\Box$

\par
\vspace{0.5cm}

Conflict of interest statement\\
The author does not have any possible conflict of interest.\\

Data availability statement\\
The manuscript has no associated data.\\




\begin{thebibliography}{99}

\bibitem{LSM} L. Aloui, S. Ibrahim and M. Khenissi, Energy decay for linear dissipative wave equations in exterior domains, J. Differential Equations {\bf 259} (2015), no. 5, 2061--2079.
\bibitem{D} W. Dan and Y. Shibata, On a local energy decay of solutions of a dissipative wave equations, Funkcial. Ekvac. {\bf 38} (1995), 545--568.
\bibitem{DA} M. Daoulatli, Energy decay rates for solutions of the wave equation with linear damping in exterior domain, Evol. Equ. Control Theory {\bf 5} (2016), no. 1, 37--59.
\bibitem{RR} R. C. Char\~ao and R. Ikehata, On the system of $2$-D elastic waves with critical space dependent damping, submitted (2025). arXiv:2503. 06854v2 [math AP] 11 Mar 2025
\bibitem{E} L. C. Evans, Partial Differential Equations, Graduate Studies in Math. {\bf 19}, AMS Providence, Rhode Island, 1998.
\bibitem{Ikawa}M. Ikawa, Hyperbolic Partial Differential Equations and Wave Phenomena; 2000. Translations of Mathematical Monographs, American Mathematical Society.
\bibitem{Ike-0} R. Ikehata, Fast decay of solutions for linear wave equations with dissipation localized near infinity in an exterior domain, J. Differ. Equ. {\bf 188} (2003), 390--405.
\bibitem{Ike-1} R. Ikehata, $L^2$-blowup estimates of the wave equation and its application to local energy decay. J. Hyperbolic Differ. Equ. {\bf 20} (2023), no. 1, 259--275.
\bibitem{IM} R. Ikehata and T. Matsuyama, $L^{2}$-behavior of solutions to the linear heat and wave equations in exterior domains, Sci. Math. Japon. {\bf 55} (2002), 33--42.
\bibitem{JJ} R. Joly and J. Royer, Energy decay and diffusion phenomenon for the asymptotically periodic damped wave equation, J. Math. Soc. Japan {\bf 70} (2018), no. 4, 1375--1418. 
\bibitem{Ma} A. Matsumura, On the asymptotic behavior of solutions of semi-linear wave equations , Publ. Res. Inst. Math. Sci., Kyoto Univ. {\bf 12} (1976), 169--189.
\bibitem{M} K. Mochizuki, Global existence and energy decay of small solutions to the Kirchhoff equation with linear dissipation localized near infinity, J. Math. Kyoto Univ. {\bf 39} (1999), no.2, 347--363.
\bibitem{MN} K. Mochizuki and H. Nakazawa, Energy decay of solutions to the wave equations with linear dissipation localized near infinity, Publ. Res. Inst. Math. Sci., Kyoto Univ. {\bf 37} (2001), 441--458. 
\bibitem{Mora} C. Morawetz, The decay of solutions of the exterior initial-boundary value problem for the wave equation, Comm. Pure Appl. Math. {\bf 14} (1961), 561--568.
\bibitem{N} M. Nakao, Energy decay for the linear and semilinear wave equations in exterior domains with some localized dissipations, Math. Z. {\bf 238} (2001), 781--797.
\bibitem{R} R. Racke, Non-homogeneous non-linear damped wave equations in unbounded domains, Math. Methods Appl. Sci. {\bf 13} (1990), 481--491.
\bibitem{SW} M. Sobajima  and Y. Wakasugi, Asymptotic expansion of solutions to the wave equation with space-dependent damping, Asymptot. Anal. {\bf 134} (2023), no. 1-2, 241--279.
\bibitem{TY} G. Todorova and B. Yordanov, Critical exponent for a nonlinear wave equation with damping, J. Differ. Equ. {\bf 174} (2001), 464--489.
\bibitem{Z} E. Zuazua, Exponential decay for the semilinear wave equation with locally distributed damping, Comm. Partial Differential Equations {\bf 15} (1990), 205--235.
\end{thebibliography}
\end{document}